\documentclass{amsart}
\usepackage{amssymb,amsmath,psfig}
\def\C{\mathbb C}
\def\H{\mathbb H}
\def\Q{\mathbb Q}
\def\R{\mathbb R}
\def\T{\mathbb T}

\def\dd{\partial}

\def\smin{\smallsetminus}

\def\h{h}
\def\o1{\overline{1}}
\def\oC{\overline{C}}
\def\oH{\overline{H}}

\def\Ga{\alpha}
\def\Gb{\beta}
\def\Ge{\varepsilon}
\def\GG{\Gamma}
\def\GD{\Delta}
\def\GL{\Lambda}
\def\Go{\omega}

\def\cJ{ J}
\def\cG{\mathcal G}
\def\cU{\mathcal U}
\def\hC{\hat{C}}
\def\tC{\widetilde{C}}
\def\tG{\widetilde{\mathcal G}}
\def\Aut{\operatorname{Aut}}

\newtheorem{thm}{Theorem}[section]
\newtheorem{lem}[thm]{Lemma}
\newtheorem{prop}[thm]{Proposition}

\theoremstyle{remark}
\newtheorem{rem}[thm]{Remark}

\def\fig#1{\psfig{figure=#1,silent=}}

\begin{document}

\title{Quantization of linear Poisson structures and degrees of maps}
\author{Michael Polyak}
\address{Department of Mathematics, The Technion, 32000 Haifa, Israel}
\email{polyak@math.technion.ac.il}
\thanks{Partially supported by the Israeli Science
Foundation grant 86/01 and the Technion grant for promotion of
research.} \keywords{Deformation quantization, Poisson structures,
Feynman graphs, configuration spaces, compactification}
\subjclass[2000]{Primary: 53D55; Secondary: 55R80, 57R35}

\begin{abstract}
Kontsevich's formula for a deformation quantization of Poisson
structures involves a Feynman series of graphs, with the weights
given by some complicated integrals (using certain pullbacks of
the standard angle form on a circle). We explain the geometric
meaning of this series as degrees of maps of some grand
configuration spaces; the associativity proof is also interpreted
in purely homological terms. An interpretation in terms of degrees
of maps shows that any other 1-form on the circle also leads to a
$*$-product and allows one to compare these products.
\end{abstract}

\maketitle

\section{Introduction}
\subsection{Kontsevich's quantization}
Recently, Kontsevich \cite{Ko} provided a deformation
quantization of the algebra of functions on an arbitrary
Poisson manifold.
The key ingredient of his construction is the quantization
formula for $\R^d$, which is then extended to the case of
general manifolds using formal geometry.
The terms of the star product are identified with certain
Feynman graphs, and their coefficients ("weights") are
given by some complicated integrals (involving certain
pullbacks of the standard volume form on a torus over the
corresponding configuration spaces).

The first non-trivial case of Kontsevich's construction is for a
linear Poisson brackets, i.e. Poisson structures on a dual of a
Lie algebra. In this paper we consider only this case (see Section
\ref{sub:problems} for a discussion).

Removal of the two external vertices splits any graph
appearing in the definition of the star product into
several connected components.
The weight of a graph is the product of weights of its
connected components.
In a linear case each of them may have at most one
cycle, which should consist of oriented edges.
It turns out that one may restrict the class of graphs
in Kontsevich's construction to trees, excluding graphs
with a non-trivial first homology.
The weights of these tree graphs coincide with the
coefficients of the usual Gutt's star product \cite{Gu}
arising from the Campbell-Baker-Hausdorff formula \cite{Ka}.
Thus Kontsevich's quantization differs from the CBH-formula
essentially only by the contributions of graphs with
one cycle, in particular of the so-called "wheels"
(see Section \ref{sub:wheels}).

\subsection{Main results and the organization of the paper}
An important problem related to Kontsevich's construction
for $R^d$ is to interpret geometrically, and show a way to
compute, the weights of graphs.
In some partial cases such computations were done in
\cite{Ka}.
However, these direct computations are complicated and
not illuminating.
One of the reasons for a lack of a geometrical
interpretation of these integrals seems to be that
some special properties of the standard angle form on
$S^1$ are significantly involved in Kontsevich's proofs.
Thus two natural questions are whether other 1-forms on
$S^1$ also lead to an associative star product, and if
yes, how do these products depend on the choice of a form.
We answer both questions in the linear case.

First we show than any 1-form on $S^1$ gives an
associative star product.
This is done by an interpretation of the associativity
as a vanishing of some cocycle in the relative top
cohomology of a grand configuration space $C_{n,3}$.
These spaces are obtained by gluing together the
configuration spaces corresponding to graphs along
their boundary strata.
These boundary strata are of two types: principal
and hidden.
The gluings along the principal strata are done
according to the Jacobi relation.
Hidden strata are zipped-up, i.e. glued to itself
by some involutions.
By choosing singular forms concentrated in one point
on a torus, one gets star products with rational
coefficients (with a straightforward estimate on the
denominators).
These coefficients can be computed by a simple
combinatorial count.

We then turn to a problem of comparison of these star products.
This is done in a similar fashion, by identifying terms in these
star products with some fixed top cohomology class of a similar
grand configuration space $C_{n,2}$. This identification involves,
however, a capping of some additional- external- boundary strata
by the configuration spaces of the wheels. Thus products for
different differential forms on $S^1$ may differ only by the
weights of the wheels. We show that the tree (i.e. 0-loop) part of
all star products coincide and are equal to Gutt's star product.
We also show that any differential form supported on a semi-circle
leads to Gutt's product.

\subsection{Related problems}
\label{sub:problems}
The point of view of this paper was heavily motivated by a similar
interpretation of the Feynman series appearing in the perturbative
Chern-Simons theory and used to produce a universal invariants of
knots and 3-manifolds, see e.g. \cite{BT, KT, Poi}. The
construction of the star product considered in this note is quite
similar (in particular, it is interesting that a similar
antisymmetry and Jacobi relations appear in the gluing
construction), although there are of course few differences (in
the Chern-Simons theory the edges are not oriented, all vertices
are trivalent- which implies a simple involution argument showing
the vanishing of the hidden strata- and the spaces after gluing
are closed). There is a clear common denominator in these
constructions; in fact, this note was initiated as a check that
the same technique may be applied. I believe that this technique
of producing grand configuration spaces by gluings is quite
powerful; it would be quite interesting to find some other
applications.

While I strongly believe that the approach used in this paper can
be applied for general Poisson structures, the zip-up argument of
Section \ref{sub:grand} for hidden faces then works only up to
degree six, so a new argument is needed to prove a conjecture
about vanishing of the weights of the hidden strata. When this
paper was already finished, the author found a modification of
this construction which works for quadratic Poisson structures (to
be discussed in a forthcoming paper). The general case, however,
still remains open.


\section{Graphs and differential operators}
\subsection{Admissible graphs}
Througout the paper we will deal with graphs $G$
such that
\begin{itemize}
\item the edges of $G$ are oriented;
\item the set of vertices of $G$ is subdivided into
      two disjoint subsets: the {\em internal} vertices
      $V^i$ and the {\em external} vertices $V^e$;
\item the set of external vertices is ordered;
\item $G$ does not contain double or looped edges.
\end{itemize}
An {\em orientation} of $G$ is an ordering of outgoing edges in
each vertex, up to a negation in any two vertices. Following
Kontsevich \cite{Ko}, we will call an oriented graph $G$ {\em
admissible} if
\begin{itemize}
\item exactly two edges start in each internal vertex,
and no edges start in any external vertex.
\end{itemize}
Since we will deal primarily with the linear case, i.e. the
Poisson brackets on a dual of a Lie algebra, we define a {\em
Lie-admissible} graph as an (oriented) admissible graph such that
\begin{itemize}
\item in each internal vertex there ends at most
      one edge.
\end{itemize}

Note that this condition implies that each connected
component of a graph obtained from $G$ by a removal
of all external vertices is either a tree, or contains
one cycle, which consists of oriented edges.

Denote by $\cG_{n,m}$ the set of all Lie-admissible graphs with
$n$ internal and $m$ external vertices. For example, $\cG_{0,2}$
is generated by the graph with two (external) vertices and no
edges. The "wedge" graph with one internal vertex connected to two
external vertices by a pair of edges (see the graph $G'$ of Figure
\ref{comp_fig}), gives two generators of $\cG_{1,2}$, which differ
by the ordering of the edges. Denote $\cG_m=\amalg_n\cG_{n,m}$.

Note that the orientation of a graph with a set $E$ of edges
induces an orientation of $\R^{|E|}=R^2\times...\R^2$ as in
\cite{Ko}: choose an arbitrary ordering of all internal vertices
and then take ordered pairs of edges starting at $1$-st, $2$-nd,
etc., vertex. For technical reasons (namely, the Jacobi relation
\eqref{Jacobi_eq} for labelled graphs) we will often need to
consider more general orderings of $E$. A {\em labelled graph} is
an admissible graph endowed with an arbitrary ordering of its
edges, compatible with the above orientation of $R^{|E|}$.
\begin{rem}\label{labels_rem}
Note that our labellings are more general than in \cite{Ko}: we
consider {\em all} orderings of edges, so for a non-oriented graph
with $n$ internal vertices and no nontrivial automorphisms the
number of labellings is $(2n)!$. Picking an orientation of $G$
cuts this number in half. In general we should also take into
account the number of automorphisms of $G$ (preserving the order
of external vertices), so the number of different labellings of an
oriented graph is $|\Aut G|^{-1}(2n)!/2$.
\end{rem}
Denote by $\tG_{n,m}$ the set of all labelled Lie-admissible
graphs with $n$ internal and $m$ external vertices, and set
$\cG_m=\amalg_n\cG_{n,m}$. There is an obvious forgetful map
$\tG_{n,m}\to\cG_{n,m}$; abusing notation, we will denote the
image of $G$ also by $G$.

\subsection{Product and composition of graphs}
There are several important operations on graphs.

\begin{figure}[htb]
\centerline{\fig{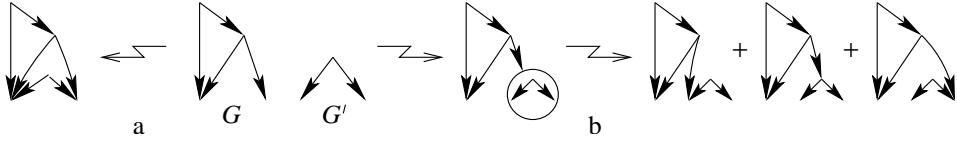,width=5.0in}}
\caption{A product $G\cdot G'$ and a composition
$G\circ_2 G'$ of graphs}
\label{comp_fig}
\end{figure}

Firstly, there is a product $\cG_{n,m}\times\cG_{n',m}
\to\cG_{n+n',m}$: we define $G\cdot G'$ to be the graph obtained
by juxtaposition of $G$ and $G'$, followed by the identification
of their external vertices. See Figure \ref{comp_fig}a. This makes
$\cG_m$ into a semigroup. In particular, $\cG_1$ is generated by
$n$-spiked "wheels", i.e. by graphs with $n$ interior vertices
forming an oriented cycle of edges and connected to the exterior
vertex by $n$ spikes, see Figure \ref{wheels_fig}a.

Secondly, there is a composition $\circ_i:\cG_{n,m}\times
\cG_{n',m'}\to\cG_{n+n',m+m'-1}$, $i=1,\dots,m$ of two admissible
graphs $G\in\cG_{n,m}$, $G'\in\cG_{n',m'}$ defined by inserting
$G'$ in $i$-th external vertex of $G$. This insertion is done in
the following way. Cut out a small disc surrounding the $i$-th
external vertex $v_i$ of $G$ and place there a small copy of $G'$.
Now, sum over all possible ways of attaching the incoming edges of
$v_i$ to the vertices of $G'$, which result in a Lie-admissible
graph. See Figure \ref{comp_fig}b. The external vertices of
$G\circ_i G'$ are ordered in a natural way:
$\{1,\dots,i-1,1',\dots,m',i+1,m\}$.

We may combine all $\circ_i$'s into one composition on the vector
space $\Q\cG$ over $\Q$ spanned by $\cG$:
\begin{equation}
G\circ G'=\sum_{i=1}^m(-1)^i G\circ_i G'.
\end{equation}

While non-associative, this operation defines a pre-Lie structure
on $\Q\cG$, i.e., for any $x,y,z\in\cG$ we have $x\circ(y\circ
z)-(x\circ y)\circ z= x\circ(z\circ y)-(x\circ z)\circ y.$ Hence
the bracket
$[G,G']=G\circ G'-G'\circ G$ defines a Lie algebra structure on
$\Q\cG$.

In a similar way, one may define a co-product on $\Q\cG$ by $\GD
G=\sum_{G'\subset G}(G/G')\otimes G'$, where the summation is over
all Lie-admissible subgraphs $G'$ of $G$ s.t. the quotient graph
$G/G'$ (obtained by shrinking $G'$ into a point) is also
Lie-admissible.

All above operations may be also defined on labelled graphs in an
obvious way.

\subsection{Differential operators assigned to graphs}
\label{diff_sec} Let $\Ga$ be a bi-vector field on an open domain
$\cU$ in $\R^d$. In the coordinates $x^i$, $i=1,2,\dots,d$ on
$R^d$ it is given by  $\Ga=\sum_{i,j}\Ga^{ij}(x)\dd_i\wedge\dd_j$,
where $\dd_i=\frac{\dd}{\dd x^i}$. To each admissible graph
$G\in\cG_{n,m}$ we assign a polydifferential operator
$B_G:(C^\infty(\cU))^{\times m}\to C^\infty(\cU)$ as in \cite{Ko}.
Namely, put $\Ga^{i(v)j(v)}$ at each internal vertex $v$ of $G$
and put the corresponding derivatives $\dd_{i(v)}$ and
$\dd_{j(v)}$ at two edges starting at $v$. Thus we have a
derivative $\dd_e$ assigned to each edge of $G$, and a tensor
$\Ga^{ij}$ at each internal vertex. We will call such an
assignment a {\em state} of $G$. Finally, put a function $f_i$ at
$i$-th external vertex, $i=1,\dots,m$. The operator $B_G$ is
defined by the following formula:

$$B_{G,\Ga}(f_1,\dots,f_m)=\sum_{\text{states}}\left(
\prod_{v\in V^i}(\prod_{e\in E(v)}\dd_e)\Ga^{i(v)j(v)}\right)
\times\left(\prod_{v\in V^e}(\prod_{e\in E(v)}\dd_e)f_v\right),$$
where $E(v)$ denotes the set of incoming edges of a vertex $v$.
For example, for the graph $G$ in the left hand side of Figure
\ref{comp_fig} we have $$B_G(f,g)=\sum_{i,j,k,l}
\Ga^{ij}\dd_j(\Ga^{kl})\dd_i\dd_k(f)\dd_l(g).$$

If $\Ga$ is linear, its second derivatives vanish, so
$B_G=0$ for any graph $G$ which contains a vertex with
at least two incoming edges.
Thus in the linear case it suffices to consider only
Lie-admissible graphs.

\subsection{Antisymmetry and Jacobi relations on graphs}
We put $G=-G'$ if two admissible graphs $G$, $G'$ differ only
by the ordering of two edges starting in some internal vertex.
Here is a graphical version of this {\em antisymmetry relation} 
($AS$):
\begin{equation}\label{AS_eq}
\vcenter{\fig{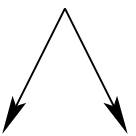}}=-\vcenter{\fig{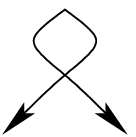}}
\end{equation}
Here and further on, in figures at each internal vertex we use
the ordering of two outgoing edges induced by the orientation
of the plane, unless specified otherwise.

We impose the following {\em Jacobi relation} on admissible
graphs:
\begin{equation}\label{Jacobi_eq}
\sum(\vcenter{\fig{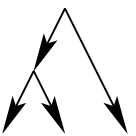}}-\vcenter{\fig{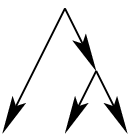}}+
\vcenter{\fig{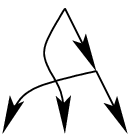}})=0
\end{equation}
In the relations above all graphs are assumed to coincide
outside the small disk shown in the figure.
There may be some incoming edges entering the disk; the
summation in the Jacobi relation is over all possible ways
to connect these incoming edges to two depicted vertices.
E.g., if there are two incoming edges, four
different ways to connect these edges to two depicted
vertices will give $4\times 3=12$ terms in the
Jacobi relation.
For Lie-admissible graphs the situation is simpler: there
may be at most one incoming edge and the unique way to
connect it, so we have
\begin{equation}\label{lJacobi_eq}
\sum(\vcenter{\fig{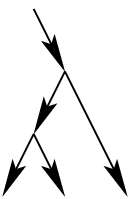}}-\vcenter{\fig{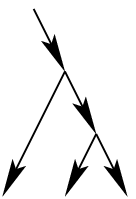}}+
\vcenter{\fig{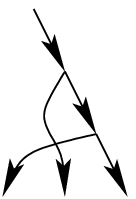}})=0
\end{equation}

There are also obvious "labelled" versions of these relations for
graphs in $\tG_{n,m}$.

Denote by $\cJ_{n,m}$ the vector space over $\Q$ generated by
graphs in $\cG_{n,m}$ modulo the $AS$ and Jacobi
relations, and set $\cJ_m=\oplus_n\cJ_{n,m}$. For $G\in\cG_{n,m}$
 denote by $[G]$ its image in $\cJ_{n,m}$ (we will also use $[G]$
 for the image of $G\in\tG_{n,m}$ under the composition of the
 forgetful map and the above quotient map).

\begin{lem}\label{Eq_lem}
Each composition $\circ_i$ of graphs factors through relations
\eqref{AS_eq}, \eqref{Jacobi_eq}, thus defines a quotient map
$\circ_i:\cJ_{n,m}\times\cJ_{n',m'}\to\cJ_{n+n',m+m'-1}$
\end{lem}
\begin{proof}
Define $[G]\circ_i[G']=[G\circ_i G']$ for any graphs
$G\in\cG_{n,m}$ and $G'\in\cG_{n',m'}$ representing the
equivalence classes $[G]$ and $[G']$ and extend it by linearity.
The fact that $[\GG\circ_i G']=0$ for any linear combination $\GG$
such that $[\GG]=0$ is immediate. It is also obvious that $[G\circ
(G'-G'')]=0$ when $G'$ and $G''$ differ by the ordering of two
outgoing edges in one vertex. The only check which may cause a
moment thought is that $[G\circ_i\GG']=0$ where $\GG'$ is a linear
combination of graphs appearing in the left hand side of equation
\eqref{Jacobi_eq}. This, however, readily follows from the fact
that in the definition of $\circ_i$ we sum over all possible
attachments of the incoming edges of $i$-th external vertex of $G$
to vertices of $\GG'$, in particular to the vertices appearing in
equation \eqref{Jacobi_eq}).
\end{proof}

The reason to call \eqref{AS_eq} and \eqref{Jacobi_eq}
"antisymmetry and Jacobi relations" is quite simple:
\begin{lem}\label{B_lem}
Let $\Ga$ be a Poisson bi-vector field. Then the map $B:G\mapsto
B_G$ factors through relations \eqref{AS_eq}, \eqref{Jacobi_eq},
thus defines a quotient map $B$ of $\cJ_m$ to polydifferential
operators.
\end{lem}
\begin{proof}
Clearly, \eqref{AS_eq} is the graphical version of the
antisymmetry relation for $\Ga$.
Equation \eqref{Jacobi_eq} with no incoming edges is exactly
a graphical version of the Jacobi identity for $\Ga$.
Similarly, \eqref{Jacobi_eq} with several incoming edges
corresponds to taking higher derivatives of the Jacobi
identity.
\end{proof}

\section{Configuration spaces and their maps}

\subsection{Configuration spaces $C_G(\oH)$ and $C_G(\R^2)$}
Consider the standard upper half-plane $\H=\{z\in\C|\Im(z)
>0\}$ as the hyperbolic plane with the Lobachevsky metric,
and the real line $\R=\{z\in\C|\Im(z)=0\}$ (together with
the point $\infty$) as the absolute.
Denote by $\oH$ the compactified space $\H\cup\R$.

For a (either labelled or not) graph $G$ with $|V^i|=n$ 
internal and $|V^e|=m$ external vertices, let $\tC_G(\oH)$  
be the space of embeddings $f:V^i\cup V^e\to\oH$
of the vertex set of $G$, such that $f(V^i)\subset\H$,
and $f(V^e)\subset\R$ preserving the order of $V^e$.
Thus we have $\tC_G(\oH)\cong\H^n\times\R^m\smin\GD$,
where $\GD$ is the union of all diagonals.
Two-dimensional group $Th(\oH)=\{z\to az+b|a\in\R_+,
b\in\R\}$ of translations and homotheties of $\oH$
acts on $\tC_G(\oH)$.
Let $C_G(\oH)$ be the quotient space $\tC_G/Th(\oH)$.
Since both $\H^n\times\R^m$ and $Th(\oH)$ are orientable,
so is $C_G(\oH)$.
This orientation does not depend on the ordering of 
internal vertices of $G$: 
if an ordering of two internal vertices is permuted, it
leads to a permutation of two copies of $\H$ in $\H^n$,
so does not change its orientation.
Further we will use the notation $C_G$ for $C_G(\oH)$
when $G$ has some external vertices, i.e. $m>0$.

For a graph $G$ with no external vertices one may consider
also a simpler configuration space: let $\tC_G(\R^2)\cong
(\R^2)^n\smin\GD$ be the space of all embeddings of the
vertex set of $G$ in $\R^2$.
Let $C_G(\R^2)$ be the quotient $\tC_G(\R^2)/Th(\R^2)$
modulo the action of a 3-dimensional group
$Th(\R^2)=\{z\to az+b|a\in\R_+, b\in\C\}$ of translations
and homotheties of $\R^2$.
Since for $m=0$ we will mostly deal with $C_G(\R^2)$ and
not with $C_G(\oH)$, in this case we shall use the notation
$C_G=C_G(\R^2)$.

\subsection{Configurations of pairs of points}
\begin{figure}[htb]
\centerline{\fig{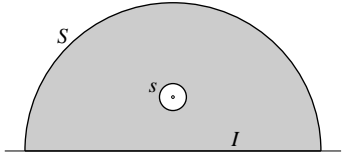,height=0.8in}}
\caption{Compactified configuration space of pairs}
\label{C2_fig}
\end{figure}

Let us describe in more details (see also \cite{Ko}) the
configuration space $C=C_e(\oH)$ of a graph $G=e$ with
two internal vertices connected by one edge.
Acting by $Th(\oH)$, we may always move the outgoing
vertex $v_1$ to $\sqrt{-1}=(0,1)$; the incoming vertex
$v_2$ then lies on $\H\smin\sqrt{-1}$.
Let us describe a compactification $\oC$ of $C$.
To compactify $C$ when $v_2$ approaches $\sqrt{-1}$, we
pick some $0<r<<1$ and cut out an open disc $|z-\sqrt{-1}
|<r$, its boundary circle $s$ encoding the direction in
which $v_2$ approaches $v_1$.
To compactify $C$ at infinity, we pick $R>>1$ and cut out
the set $|z|> R$, its boundary semicircle $S$ encoding the
direction in which the second point approaches infinity.
Finally, to compactify $C$ when $v_2$ approaches $\R$,
we add an interval $I=\{z\in\C|-R<\Re(z)\le R,\Im(z)=0\}$.
See Figure \ref{C2_fig}.
Notice also, that the space $C_e(\R^2)\cong S^1$, and
may be identified with the circle $s$ after an
identification of $\R^2$ with $T_{\sqrt{-1}}\H$.

Suppose now that one of the vertices of a one-edge
graph is external.
The corresponding configuration space is an interval;
in the model above it may be identified either with
the semicircle $S$ (if $v_1$ is external) or with the
interval $I$ (if $v_2$ is external).
Finally, suppose that both vertices are external.
In this case the space $C_e(\oH)$ is a point and may be
identified in the model above with one of the endpoints
of the interval $I$, depending on the direction of the
edge (or, equivalently, with the endpoints of the
semicircle $S$).

Thus various boundary strata of $\oC$ actually correspond
to all possible configuration spaces of pairs of points.

\subsection{Angle maps}
\begin{figure}[htb]
\centerline{\fig{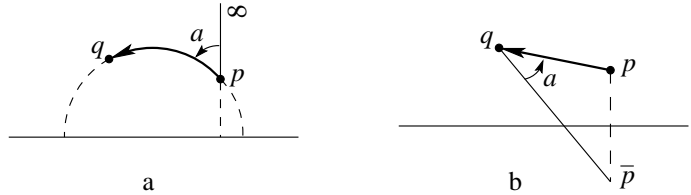,height=1in}}
\caption{Angle maps}
\label{angle_fig}
\end{figure}

The space $\oC$ has the homotopy type of $S^1$.
We will call a map $a:\oC\to S^1$ an {\em angle map},
if it is identity on $s$, maps the whole of $S$ to $0$,
and maps $I$ onto $S^1$ surgectively.
An angle map used in \cite{Ko} is defined as follows.
Let $(p,q)\in\oH^2$ be a pair of points.
Compute the angle between the direction from $p$ to
$\infty$ and the direction from $p$ to $q$ in the
hyperbolic geometry on $\oH$ (see Figure \ref{angle_fig}a).
A simple check shows that this map is $Th(\oH)$-invariant
and indeed determines an angle map.

Here is another simple example of an angle map: for
$(p,q)$ as above, take the angle between the direction
from $q$ to $\bar{p}$ and from $q$ to $p$ in the
Euclidean geometry on $\C$ (see Figure \ref{angle_fig}b).

We will fix some angle map $a$ for the remaining part of
the paper.

\subsection{Gauss maps and compactifications of
configuration spaces}
\label{sub:gauss}
For an edge $e$ of a labelled graph $G$, an evaluation
map $ev_e:C_G\to\oC$ is obtained by forgetting all
vertices of $G$ except for the ends of $e$.
A Gauss map $\phi_e:C_G\to S^1$ is a composition
of $ev_e$ with the angle map $a$.
Taking the map $\phi_e$ for each edge $e_1,\dots,e_k$
and using the ordering of edges, we get the product
Gauss map $\Phi_G=\prod_e\phi_e:C_G\to\T^k$.
We shall construct compactification $\oC_G$ of $C_G$ so
that $\Phi_G$ will extend to $\oC_G$.
The minimal construction is to take the closure of the
graph of $\phi_G$ in $(\H^n\times\R^m)\times\T^k$;
unfortunately it is rather hard to describe explicitly.
Another extreme ("overblown") solution is to blow up
all diagonals of $\H^n\times\R^m$; this compactification
was used in \cite{Ko}.
The map $\Phi_G$ obviously extends to any of the above
compactifications.

Following \cite{KT}, we will consider an intermediate solution: to
blow up only the diagonals of $\H^n\times\R^m$ corresponding to
subgraphs of $G$ which are 2-connected (this is somewhat smaller
than in \cite{Poi}, where blow ups were done for all connected
subgraphs). Here a graph is {\em 2-connected} (or rather
vertex-2-connected), if it remains connected after a removal of
any one of its vertices. In particular an edge is 2-connected. The
resulting space $\oC_G$ is a manifold with corners. The boundary
strata look as follows. Each codimension one face of $\oC_G$
corresponds to a collision of all vertices of a
vertex-one-connected subgraph $\GG$ of $G$. Both the face and the
restriction of the Gauss map to it have a product structure
$\Phi_{G/\GG}\times\Phi_\GG:C_{G/\GG}\times C_\GG\to
\T^{l-k}\times\T^k$, where $l$ and $k$ is the number of edges of
$G$ and $\GG$ respectively. Here the second factor encodes the way
the vertices of $\GG$ approach each other, seen through a
"magnifying glass" (called a {\em screen} by Fulton and
MacPherson). More generally, a stratum of codimension $k$
corresponds to a collection of $k$ vertex-one-connected subgraphs
$\GG_i$ of $G$, each two of which are either nested (i.e.
contained in one another) or disjoint. Such a stratum is adjacent
to all strata of codimension $k-1$ obtained by removing one of
$\GG_i$'s from the collection.


\subsection{Constructing the grand configuration space}
\label{sub:grand} Consider the disjoint union
$\tC_{n,m}=\amalg_{G\in \tG_{n,m}}\oC_G$ of configuration spaces
and the map $\Phi_{n,m}=\amalg_G\Phi_G:\tC_{n,m}\to\T^{2n}$. The
space $\tC_{n,m}$ has no top degree cohomology since each
component $\oC_G$ has faces. We want to get rid of the faces of
$\tC_{n,m}$ by all means: gluing them together, capping-off, or
relativizing, to obtain a space $C_{n,m}$ with a much smaller
boundary, but still with a well-defined map
$\Phi_{n,m}:C_{n,m}\to\T^{2n}$.

Some boundary strata are easy to take care of; namely,
we shall not care for any strata which are
\begin{itemize}
\item of codimension two or more in $\tC_{n,m}$, or
\item mapped to $\T^{2n}$ with a loss of dimension
      (i.e. the difference between the dimensions
      of the source and the image) two or more, or
\item mapped to the union $O$ of all coordinate
      hyperspaces $\{\Phi_i=0\}\subset\T^{2n}$.
\end{itemize}

We will simply relativize these strata including them
in a singular locus $D_{n,m}$ and considering then
$\Phi_{n,m}:(C_{n,m},D_{n,m})\to (\T^{2n},O)$.
Let us describe the gluing/relativizing process for
various types of boundary faces of $\tC_{n,m}$.

\begin{figure}[htb]
\centerline{\fig{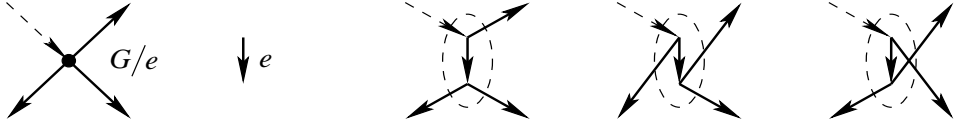,width=5.0in}}
\caption{Gluing principal faces together}
\label{main_fig}
\end{figure}

Firstly, there are the {\em principal internal faces},
where two internal vertices of a graph $G$ connected by
one edge $e$ collide together.
Note that the same principal face $C_{G/e}\times S^1$
appears in each of the three graphs related by the Jacobi
relation, as shown in Figure \ref{main_fig}, where the
colliding vertices are indicated by a dashed circle (and
a possible incoming edge is also dashed).
Thus we can glue these faces of three graphs of Figure
\ref{main_fig} together.

Secondly, there are the {\em hidden internal faces}
$C_{G/\GG}\times C_\GG$, where several internal vertices of a
connected subgraph $\GG$ of $G$ (with more than one edge) collide
together. A simple dimension counting shows that unless the number
of edges starting in $\GG$ and ending in $G\smin\GG$ is two or
three, the loss of dimension under the map is at least two on one
of the two factors, so the face is readily relativized. Also,
unless $\GG$ is 2-connected, $C_\GG$ loses at least one dimension
under $\Phi_G$ (since the dilatations of one of the parts into
which a cut-vertex splits $\GG$ preserve $\Phi_G$), so such a face
can again be relativized. The remaining hidden faces\footnote{This
is actually the only place where we have to use Lie-admissibility
of the initial graph, and thus the linearity of the Poisson
structure} have an orientation reversing involution which
preserves $\Phi_G$, thus can be zipped-up (i.e. glued to itself by
this involution) following an idea of \cite{Ko2}. Indeed, pick a
vertex $v$ in $\GG$ with an edge ending in $G\smin\GG$. Any vertex
of a Lie-admissible graph has valence at most three; moreover,
$\GG$ is 2-connected, thus any vertex $v\in\GG$ with an edge
ending in $G\smin\GG$ has precisely two edges $e_1=(v_1,v)$ and
$e_2=(v,v_2)$ with $v_1,v_2\in\GG$. The involution of $v$ w.r.t.
the middle of the interval connecting $v_1$ and $v_2$ (i.e. $v\to
v_1+v_2-v$, if we think about $v$, $v_1$ and $v_2$ as points in
$\R^2$), followed by the reordering of $e_1$ and $e_2$ reverses
the orientation of this boundary face, but preserves $\Phi_G$.

Thirdly, there are the {\em infinity faces}, where some
internal vertices of a connected subgraph $\GG$ of $G$
approach the absolute line $\R$, but remain distant from
the external vertices.
Then any edge starting in $\GG$ and ending in $G\smin\GG$
is mapped to 0, thus such faces can again be relativized.

Finally, there are the {\em external faces}, where the
vertices of a connected subgraph $\GG$ of $G$ which
collide together contain $k>0$ external vertices.
Again, if there is an edge starting in $\GG$ and ending
in $G\smin\GG$, it is mapped to 0, thus such faces can
be relativized.
In the remaining cases the situation depends on $k$.
If $k>2$ or $k=1$, then the map $\Phi_{G/\GG}$ or
$\Phi_\GG$ respectively has a loss of at least two
dimensions, so again can be relativized.
We are left, however, with the case $k=2$.
We can still get rid of at least some of these faces.
Namely, if the graph $G/\GG$ has a double edge $e_1=e_2$,
we zip up a pair of faces corresponding to (a collapse of
a subgraph $\GG$ of) two graphs $G$ and $G^*$ which differ
only by the ordering of $e_1$ and $e_2$.
Indeed, the maps $\Phi_G$ and $\Phi_{G^*}$ differ only by
a transposition of two factors of $S^1$ in $\T^{2n}$, so
coincide up to a sign on the diagonal $e_1=e_2$.

To sum it up: the only boundary faces of $C_{n,m}$ (which are not
included in the singular locus $D_{n,m}$) are the external faces,
corresponding to subgraphs $\GG\in\cG_{k,2}$ of a graph
$G\in\cG_{n,m}$ which contain exactly two external and any number
$k$ of internal vertices, all edges starting in $\GG$ also end
there, and $G/\GG$ does not have double edges. In other words,
\begin{prop}
\label{prop:dd}[cf. Kontsevich \cite{Ko}] The boundary faces of
$C_{n,m}$ restricted to $\oC_G$ correspond to pairs
$\GG'\in\tG_{n-k,m-1}$ and $\GG''\in\tG_{k,2}$ such that
$G=\GG'\circ_i\GG''$ for some $i$. The restriction of $\Phi_{n,m}$
to such a face is
\begin{equation}
\label{eq:dd}
\Phi_{n-k,m-1}\times\Phi_{k,2}:C_{\GG'}\times
 C_{\GG''}\to\T^{2(n-k)}\times\T^{2k}
\end{equation}
\end{prop}
We will use this crucial fact in Section \ref{sub:3v}
below, and then again in Section \ref{sub:compare-bad}.

\subsection{Example: wedge graphs}
\label{sub:3pt} 
\begin{figure}[htb]
\centerline{\fig{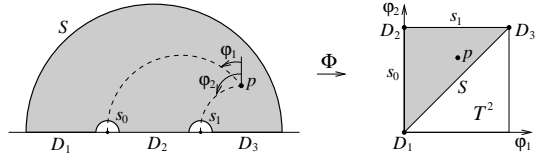,height=0.8in}}
\caption{Compactified space of the wedge graph and its map to $\T^2$}
\label{fig:3pt}
\end{figure}

There are two "wedge" graphs $G$ in $\cG_{1,2}$
which differ only by the ordering of the edges; denote by $G$ the
one with the natural order of the edges, and by $G^*$ the other.
To visualize the space $C_G$, it is convenient to fix two external
vertices in $0,1\in\R$ by the action of $Th(\oH)$. The remaining
internal vertex is free to run in the complement of these two
points in $\H$, or, equivalently, in the complement of
$0,1,\infty$ on $\oH$, so $C_G\cong\oH\smin\{0,1,\infty\}$. The
compactification $\oC_G$ is shown on Figure \ref{fig:3pt}.

The small semi-circles $s_0$, $s_1$ around $0$, $1$ correspond to
the faces where the internal vertex approaches one of the external
vertices; $\Phi_G(s_0) =S^1\times\{0\}$ and $\Phi(s_1)=\{0\}\times
S^1$. The singular locus $D$ is the union of the pieces $D_1$,
$D_2$, $D_3$ of the absolute line $\R$ together with $s_0$
and $s_1$; $\Phi_G(D)\subset O$ (recall that $O\subset\T^2$ is
the union of the coordinate hyperspaces, i.e.
$O=\{\Phi_1=0\}\cup\{\Phi_2=0\}$).

The large semicircle $S$ corresponds to the infinite face where
two external vertices collide together (or, equivalently, the
internal vertex approaches the infinity on $\H$); $\Phi_G(S)$ is
the diagonal $\GD\subset\T^2$. This is the only non-relativized
boundary face of $\oC_G$.

The space $\oC_G$ is a hexagon with the boundary edges
$D_1$, $s_0$, $D_2$, $s_1$, $D_3$ and $S$; the relativised
space $(\oC_G,D)$ is a disc with the boundary $S$.
See Figure \ref{fig:3pt}.

The space $\oC_{G^*}$ of the graph $G^*$
is the same; in order to keep track of the spaces we will denote
the corresponding locus by $D^*$ and the boundary face by $S^*$.

Each of the open cells $C_G$ and $C_{G^*}$ has its
unique 2-dimensional cohomology class, and thus the
top cohomology of $C_G\cup C_{G^*}$ (or, equivalently,
$H^2(\oC_G\cup\oC_{G^*},D\cup D^*\cup S\cup S^*)$)
is 2-dimensional with two generators, which we may
identify with $G$ and $G^*$.

\begin{figure}[htb]
\centerline{\fig{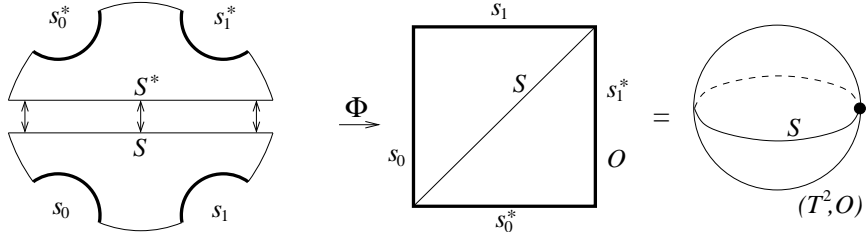,width=4.5in}}
\caption{Gluing $\oC_G$ and $\oC_{G^*}$ together and
mapping them to $S^2$}
\label{fig:wedge}
\end{figure}

The map $\Phi_{G^*}$ differs from $\Phi_G$ by the transposition of
two copies of $S^1$ in $\T^2$. In particular, on the diagonal
$\GD$ of $\T^2$ both maps coincide up to a sign. Thus we may glue
$\oC_G$ to $\oC_{G^*}$ along the remaining (non-relativised)
boundary faces $S$ and $S^*$. The gluing however changes the
cohomology: we have to add the correponding relation $G^*=-G$,
making $H^2(C_{1,2},D)$ into 1-dimensional. Indeed,
$(C_{2,1},D_{2,1})$ is simply $S^2$, obtained by gluing two discs
$(\oC_G,D)$ and $(\oC_{G^*},D^*)$ along their boundaries $S$ and
$S^*$. Note, that it is mapped by $\Phi$ also to $S^2$, since
contracting $O\subset\T^2$ into a point we get $(\T^2,O)\cong
S^2$. See Figure \ref{fig:wedge}.

This simple example prepares us for a general
calculation of the top cohomology of the grand
configuration space below.

\subsection{Cohomology of the configuration space}
\label{ssec:cohom} Let us compute the top cohomology of
$\oC_{n,m}$. Each cell $C_G$ has a unique top cohomology class-
the fundamental class. Thus
$H^{2n+m-2}(\tC_{n,m},X_{n,m})=H^{2n+m-2} (\amalg_G
C_G)=\oplus_{G\in\tG_{n,m}} \Q G$, where $X_{n,m}$ is the union of
all faces of $\tC_{n,m}$. The space $C_{n,m}$ was obtained from
$\tC_{n,m}$ by some gluings. Each time we glue some top
dimensional cells of $\tC_{n,m}$ along a face, it adds a relation
on the generators of the cohomology. The gluings along the primary
internal faces were made according to the Jacobi relation. The
zip-ups of the hidden strata were made according to the
antisymmetry relation. No other relations were introduced during
these zip-ups, since we glued together several graphs which are
identical under $\Phi_{n,m}$ up to a sign; half of them are
negated by the map $\Phi_{n,m}$, so the total sum is 0. We used
the antisymmetry relation again to zip up the external faces
corresponding to graphs with double edges. Thus
$H^{2n+m-2}(C_{n,m},\dd C_{n,m}\cup D_{n,m})$ is the space of
labelled Lie-admissible graphs modulo the Jacobi and the
antisymmetry relations. This is not quite $\cJ_{n,m}$, since the
graphs are labelled, but not too different: forgetting a labelling
gives a surjective map to $\cJ_{n,m}$.

Thus a class in $H^{2n+m-2}(C_{n,m},\dd C_{n,m}\cup D_{n,m})$ (or
rather its image $Z_{n,m}\in\cJ_{n,m}$ under the forgetful map) of
a top-degree form $w$ can be written as
\begin{equation}
\label{eq:Z_nm} Z_{n,m}=\sum_{G\in\tG_{n,m}} W_G [G]
\end{equation}
where $[G]$ is the equivalence class of $G$ in $\cJ_{n,m}$, and
the {\em weights} $W_G$ can be calculated as $W_G=\int_{C_G} w$.
For what follows it will be important to note that if $w$ is
chosen so that $w|_{C_G}=d\Gb$, then by the Stokes theorem we can
compute the weights also as $W_G=\int_{\dd C_G}\Gb$.

\begin{rem}
There is another version of this construction in which one works
with unlabelled graphs, which are mapped not to $\T^{2n}$, but to
its quotient under the action of $S_{2n}$ transposing the copies
of $S^1$. While it is somewhat cleaner and requires less auxiliary
data, it is harder to visualize topologically. There is, however,
a simple combinatorial counterpart, which allows one to reduce the
summation in \eqref{eq:Z_nm} to unlabelled- and even further to
connected- graphs, as will be explained in Section \ref{sub:two}
below.
\end{rem}

Further on we will be mostly interested in $m=2$ and $m=3$. For
these values of $m$ we have $\dim(\tC_{n,2})=2n$, and
$\dim(\tC_{n,3})=2n+1$, so $\tC_{n,2}$ is mapped by $\Phi_{n,2}$
to the torus of the same dimension, while $\tC_{n,3}$ has one
excessive dimension.

\subsection{Two external vertices}\label{sub:two}
For $m=2$ we fix a non-trivial 1-form $\Go$ on $S^1$ normalized by
$\int_{S^1}\Go=1$, and take for $w$ the pull-back
$\Phi_{n,m}^*(\Go^{2n})$ of the volume form
$\Go^{2n}\in\GL^{2n}(\T^{2n})$. So we set
$$W_G=\int_{C_G} \Phi_{n,m}^*(\Go^{2n})$$
in \eqref{eq:Z_nm}. Since different labellings of a graph
$G\in\cG_{n,2}$ lead to even transpositions of copies of $S^1$ in
$\T^{2n}$ and $\Go^{2n}$ is invariant under such transpositions,
the weight $W_G$ is the same for all labellings of $G$. Thus,
instead of summing the weights over graphs in $\tG_{n,2}$, we may
sum over $G\in\cG_{n,2}/AS$ (picking any orientation and any 
labelling compatible with
 it to define $\Phi_G$ and $W_G$), counting 
each $G$
 with a multiplicity being the number of its different labellings,
i.e. $(2n)!|\Aut G|$, see Remark \ref{labels_rem} (the additional factor 
of 2 disappeared since instead of counting each graph in $\cG_{n,2}$ 
twice with both orientations we take a quotent by the AS relation).

Summing over all $n$ with a formal parameter $h^n$ and an
appropriate normalization factor, we denote by
 $Z$ the 
following generating function:
\begin{equation}
\label{eq:Z} Z=\sum_n \frac{h^nZ_{n,2}}{(2n)!}=
\sum_
{G\in\tG_2}\frac{h^{|G|} W_G}{(2n)!} [G]=
\sum_{G\in\cG_2/AS}\frac{h^{|G|} W_G}{|\Aut G|} [G]
\end{equation}
where $|G|$ is the number of the edges of $G$. Noticing that
$\Phi_{G'\cdot G''}=\Phi_{G'}\times \Phi_{G''}$ for any
$G',G''\in\cG_2$, we get
\begin{lem}
The weights are multiplicative, i.e.
$W_{G'\cdot G''}=W_{G'}W_{G''}$ for any $G',G''\in\cG_2$.
\end{lem}

There are still quite many terms in each degree in the definition
of $Z$, since we sum over all graphs in $\cG_2$. There is a
standard trick of taking the logarithm, allowing one to count only
connected graphs. Denote by $\cG_2^c$ the subset of graphs in
$\cG_2$ which remain connected after a removal of both external
vertices.

\begin{prop}
\label{prop:ln}
$$\ln(Z)=\sum_{G\in\cG_2^c/AS}
\frac{h^{|G|} W_G}{|\Aut G|}[G]$$
\end{prop}
Indeed, for distinct graphs $G_1,\dots,G_k$ we have
$$|\Aut (G_1^{i_1}\cdot \dots G_k^{i_k})|=i_1!\dots
i_k!\ |\Aut G_1|^{i_1}\dots |\Aut G_k|^{i_k}$$
which together with the multiplicativity of the weights
proves the lemma.

\subsection{Coefficients of $Z$ in small degrees}
\label{sub:small} Let us check the coefficients of $Z$ for small
values of $n$. The only graph in $\cG_{0,2}$ has no internal
vertices (and no edges); its weight is 1. For $n=1$ we have two
wedge graphs $G$, $G^*$, as discussed in Section \ref{sub:3pt}.
Since $|\Aut G|=|\Aut G^*|=1$, we get $1/2(W_G-W_{G^*})[G]$ for
the coefficient of $h$ in $Z$. How can we find $W_G-W_{G^*}$
without direct calculations? Note that the map
$\Phi_{1,2}:(C_{1,2},D_{1,2})\to(\T^2,O)$ has a well-defined
relative degree. It can be computed by pulling back a volume form
$\Go^2$ on $\T^2$ and integrating it over $C_G\cup C_{G^*}$, so
$\deg(\Phi_G)=\int_{C_{1,2}}\Phi_G^*(\Go^2)= W_G-W_{G^*}$. Another
way to compute $\deg(\Phi_G)$ is to count the algebraic number of
preimages $\Phi_G^{-1}(r)$ of any regular value $r$. Choosing some
reasonable regular value, we see that there is only one preimage
(e.g., for $r=(\pi/2,\pi)\in\T^2$ we get $\Phi_G^{-1}(r)=1+i\in
C_G$); thus the coefficient of $h$ in $Z$ does not depend on the
choice of $\Go$ and is $1/2$.

We will use similar arguments in Section \ref{sub:compare-good}
below.

\subsection{Three external vertices}
\label{sub:3v} For $m=3$ the pull-back
$\Gb=\Phi_{n,m}^*(\Go^{2n})$ of a volume form
$\Go^{2n}\in\GL^{2n}(\T^{2n})$ to $C_{n,3}$ is a form of degree
$2n$, while $\dim(C_{n,3})=2n+1$. Note that
$w=d\Gb=d\left(\Phi_{n,m}^*(\Go^{2n})\right)=
\Phi_{n,m}^*(d(\Go^{2n}))=0$. Thus $w$ represents the zero
cohomology class in the top cohomology $H^{2n+1}(C_{n,3},\dd
C_{n,3}\cup D)$. Let us now compute the same zero differently-
this will give a non-trivial equation on $Z$. By Proposition
\ref{prop:dd} for $m=3$, each boundary face of $C_{n,3}$
corresponds to $\GG'\circ\GG''$ with $\GG'\in\cG_{n-k,2}$,
$\GG''\in\cG_{k,2}$, and the restriction of $\Phi_{n,3}$ to it has
a product structure
$$\Phi_{\GG'}\times\Phi_{\GG''}:C_G\times C_\GG\to
\T^{2(n-k)}\times\T^{2k}$$
Therefore, as in Section \ref{ssec:cohom} above, the image
$0=w=d\Gb$ in $H^{2n+1}(C_{n,3},\dd C_{n,3}\cup D_{n,3})$ of
$\Gb\in H^{2n}(\dd C_{n,3},D)$ may be calculated as
\begin{multline*}
\sum_k\sum_{\GG'\in\cG_{n-k,2}}\sum_{\GG''\in\cG_{k,2}}
\frac1{|\Aut(\GG'\circ\GG'')|}\int_{C_{\GG'}\times
C_{\GG''}}\Phi_{n-k,2}^*(\Go^{2(n-k)})\times
\Phi_{k,2}^*(\Go^{2k})\ [\GG'\circ\GG'']\\ =\sum_k
\left(\sum_{\GG'\in\cG_{n-k,2}}\frac{W_{\GG'}}
{|\Aut\GG'|}[\GG']\right)\circ\left(\sum_{\GG''\in
\cG_{k,2}}\frac{W_{\GG''}}{|\Aut\GG''|}[\GG'']\right)
\end{multline*}
and we conclude that
\begin{prop}
\label{prop:ZZ}
$Z\circ Z=0$.
\end{prop}

\section{Quantization of Poisson structures}
\label{sec:quant}

\subsection{Constructing star products}
Let $A=\GG(M,\mathcal{O}_M)$ be the algebra of
$C^\infty$-functions on a manifold $M$.
A {\em star product} is an associative $\R[[\h]]$-linear
product on $A[[\h]]$ which, for functions
$f,g\in A\subset A[[\h]]$ is a deformation of the standard
product $fg$:
\begin{equation}\label{star_eq}
f*g=fg+B_1(f,g)\h+B_2(f,g)\h^2+\dots,
\end{equation}
where $B_i:A\times A\to A$ are bidifferential operators.
The product then extends to $A[[h]]$ by linearity.
A {\em deformation quantization of a Poisson structure}
on $M$ is a construction of a star product, first term
$B_1$ of which coincides with the (half of) Poisson bracket
$\{\ ,\ \}:A\times A\to A$.

Kontsevich \cite{Ko} have shown that such a quantization
always exists.
His approach is based on an explicit construction of a
star product on (an open domain of) $\R^d$, which then
extends to an arbitrary $M$.

The dual of a Lie algebra has a natural Poisson structure
which is linear, i.e. the components $\Ga^{ij}$ of the
Poisson tensor $\Ga=\Ga^{ij}\dd_i\wedge\dd_j$ depend
linearly on the coordinates $x^i$.
Below we provide a modification of Kontsevich's
original construction for linear Poisson structures.

Let $\Ga$ be a linear Poisson bi-vector field on an open
domain $\cU\subset\R^d$, defining a Poisson structure
on the algebra $A$ of smooth functions on $\cU$.
Fix an arbitrary 1-form $\Go$ on $S^1$ and define
a star product \eqref{star_eq} using the series $Z$
of equation \ref{eq:Z} and the map $B$ of Section
\ref{diff_sec} from graphs to bidifferential operators
by $$f*g=B_{Z,\Ga}(f,g)$$ for $f,g\in A$.

Returning to the computations of the coefficients of
$Z$ in Section \ref{sub:small} we see that
$f*g=f\cdot g+\frac{h}{2}\{f,g\}+\dots$
Moreover, translating the statement $Z\circ Z=0$ of
Proposition \ref{prop:ZZ} to the language of
differential operators, we readily obtain

\begin{thm}
A star product $f*g=B_{Z,\Ga}(f,g)$ is associative.
\end{thm}

The above star product {\em a priori} depends on the initial
choice of the 1-form $\Go$ on $S^1$. It is easy to see directly
from the definition that Kontsevich's product \cite{Ko}
corresponds to $\Go=d\phi$ being the standard angle form on $S^1$.

What can be said about other 1-forms? This is the subject of the
remaining part of this note.

\subsection{Digression on degrees of maps}
\label{sub:compare-good}
If $C_{n,2}$ would have no boundary
faces apart from the singular locus $D_{n,2}$, then (similarly to
the case $n=1$ in Section \ref{sub:small}) $Z_{n,2}$ could be
interpreted as the relative degree of a map
$\Phi_{n,2}:(C_{n,2},D_{n,2})\to(\T^{2n},O)$.

Let me clarify this point: while usually one thinks about the
degree $\deg f$ of a map $f:M\to N$ between two manifolds being a
number, one can equivalently think that it is a vector $\deg f
[M]$ in a 1-dimensional vector space (the top cohomology of $M$,
generated by the fundamental class $[M]$ of $M$). Equivalently,
$\deg f$ can be defined as a pull-back $f^*([N])$ of the
fundamental class of $N$, since $f^*([N])=\deg f [M]$. In our
setting $M$ is not a manifold, but a stratified space, so its top
cohomology is not 1-dimensional, but is generated by graphs. While
one can not any more define the degree as a number, the definition
of degree as $f^*([N])$ generalizes to this situation, with all
usual advantages of presenting a topological value as a degree.

In particular, we would be able to make two conclusions:
\begin{itemize}
\item All 1-forms $\Go$ give the same star product;
\item The coefficients of this star product are
      rational (with the denominators arising from
      the normalization factors $1/(2n)!$).
\end{itemize}

Unfortunately, this is not the case: the above degree is not
well-defined, since in general $C_{n,2}$ has some boundary faces
(so when we change $\Go$, the value of $Z$ may change due to a
"flow" through these faces). However, at least some partial
steps on this way may be achieved.

\subsection{Comparison of star products: 0-loop part}
\label{sub:compare-bad}

To understand the boundary faces, note that there is an additional
grading on the configuration space by the loop number of graphs.
Here the {\em loop number} $\ell(G)$ of a graph $G$ is a number of
connected components appearing after a removal from $G$ of both
external vertices. Indeed, the gluing relations \eqref{AS_eq} and
\eqref{Jacobi_eq} preserve the loop number. Thus $C_{n,2}$
consists of non-connected components $C_{n,2}^\ell$ corresponding
to graphs with $\ell(G)=\ell$. Accordingly, any star product
splits into pieces enumerated by the loop number of the
corresponding graphs. Gutt's star product corresponds to the
0-loop piece:
\begin{prop}
The 0-loop part of a star product do not depend on the choice of a
1-form $\Go$ on $S^1$ and coincides with the standard Gutt
product.
\end{prop}
\begin{proof}
The first statement follows from the fact that the 0-loop part
$C_{n,2}^0$ of the configuration space has no non-relativized
boundary faces, so the scheme discussed in Section
\ref{sub:compare-good} above works and the 0-loop part $Z^0_{n,2}$
of $Z_{n,2}$ is, up to a normalization by $1/(2n)!$, a relative
degree of the restriction of $\Phi_{n,2}$ to
$(C^0_{n,2},D^0_{n,2})\to(\T^{2n},O)$. Thus the 0-loop part of all
star products are equal. Now results of Kathotia \cite{Ka}, who
verified that for the standard angle 1-form on $S^1$ the
coefficients of tree graphs coincide with the Gutt product, imply
the theorem.

One can also verify this directly, in a purely combinatorial
manner, computing the coefficients by counting the algebraic
number of preimages in $C^0_{n,2}$ of any regular value
$r\in\T^{2n}$ (so up to normalization factor $1/(2n)!$ we will get
integral coefficients). Note, that this counting is to be
performed with labelled graphs; indeed, to pass to unlabelled
graphs in Section \ref{sub:two} we used symmetries, thus would
have to deal with (not regular) values on the diagonal of
$\T^{2n}$. A lengthy, but rather routine calculation (using
Proposition \ref{prop:ln} for $\ln(Z)$) shows that it indeed gives
Gutt's product.
This is illustrated in the Example \ref{sub:preimage} below.
\end{proof}

\subsection{Example}\label{sub:preimage}
\begin{figure}[htb]
\centerline{\fig{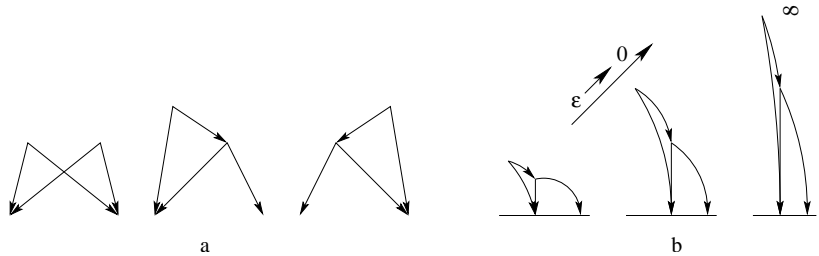,height=1.3in}}
\caption{Graphs in $\cG_{2,2}$ and preimages of points approaching
the diagonal} \label{fig:n2}
\end{figure}

Let us compute the coefficients of 0-loop graphs for $n=2$, see
Figure \ref{fig:n2}a. The easiest way to perform calculations is
to pick a regular value $r\in\T^4$ very close to the diagonal.
Say, let us take $r=(\pi,\pi+\Ge,\pi+2\Ge,\pi+3\Ge)$. When
$\Ge\to0$, the value approaches the diagonal and internal vertices
of graphs in the preimages tend to the infinity on $\H^2$. For
small $\Ge$ the edges connecting these vertices to the external
vertices are almost Euclidean straight lines, as illustrated in
Figure \ref{fig:n2}b. Thus calculations may be done for simplicity
in the Euclidean geometry.

The first graph of Figure \ref{fig:n2}a is the square of the wedge
graph (it splits into two copies of the wedge graph after the
removal of the external vertices), thus its coefficient should be
$\frac1{2!}(\frac12)^2=\frac18$ in view of Proposition \ref{prop:ln}.
This is equally easy to verify diectly: out of the total of 12
different labellings of 4 edges by $\{1,2,3,4\}$ exactly 3
(namely $(e_1,e_2,e'_1,e'_2)=(1,2,3,4)$, $(1,3,2,4)$, $(1,4,2,3)$
and $(2,3,1,4)$) give a preimage. It is also easy to check that
for all of them the local degree is $+1$, so the coefficient is
indeed $\frac1{4!}3=\frac18$.

\begin{figure}[htb]
\centerline{\fig{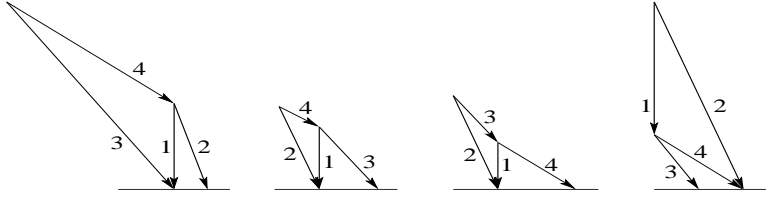,height=1in,width=4in}}
\caption{Labellings corresponding to preimages}
\label{fig:preimage}
\end{figure}

As for the second graph of Figure \ref{fig:n2}a, out of the total
of 24 different labellings exactly 4 give a preimage, see Figure
\ref{fig:preimage}.
This is due to elementary geometrical/combinatorial reasons.
Consider the internal vertex $v$ for which one of the starting
edges ends in the second internal vertex. Obviously, both edges
starting in $v$ should eminate at angles which are either both
smaller or both greater than the angle of the edge $e$ connecting
their ends. In other words, this pair of edges should have labels
either both smaller or both greater than the label of $e$; see
Figure \ref{fig:preimage}. It is easy to check that there are only
4 such labellings, depicted in Figure \ref{fig:preimage}. It is
also easy to check
 that the local degree of exactly one of them
(the last labelling in Figure \ref{fig:preimage}) is $-1$,
thus the coefficient is $\frac1{4!}(3-1)=\frac1{12}$.

By a similar inductive combinatorial arguments one can calculate
the numbers of
 preimages of 0-loop graphs with higher $n$ (getting
modified Bernulli
 numbers). Direct visualization of preimages
becomes hard for
 $n>2$ since numbers of labellings (as well as
denominators) grow
 too fast: for $n=4$ it is already $8!=40320$.

\subsection{Graphs with loops}
\label{sub:loops}

It remains to study the contribution of graphs with non-zero loop
numbers. For some choices of an angle form it is easy:

\begin{prop}
Let $\Go$ be any 1-form on $S^1$ supported on a semicircle
$\{z=exp(i\psi),a<\psi<a+\pi\}$. Then the corresponding star
product is the Gutt product.
\end{prop}
\begin{proof}
Without a loss of generality, suppose that the form is supported
on the lower semicircle $\{|z|=1, \Im(z)< 0\}$. Now, the weight of
any graph which contains an oriented cycle is zero, since such a
cycle should contain an edge directed upwards. But in the linear
case each graph with non-zero loop number contains such a cycle.
Thus only 0-loop graphs contribute to the star product.
\end{proof}

In other cases, the situation is more complicated and at present I
fall short of computing the weights of graphs with non-zero loop
numbers for arbitrary 1-forms on $S^1$. The remaining part of this
note thus has a preliminary character. Its purpose is to explain
how a computation of weights of graphs with loops (and two
external vertices) can be done instead using so-called wheels,
which are graphs of the type depicted in Figure \ref{wheels_fig}c.
This gives a possibility to compute Kontsevich's weights of
graphs with loops via a purely combinatorial count (see Section 
\ref{sub:twospike} below for an example).

\subsection{Wheels}
\label{sub:wheels}
The above reasoning may be repeated for general $n$.
Indeed, let us return to the study of the boundary faces of $C_{n,2}$.
By Proposition \ref{prop:dd}, these faces correspond to
decompositions $G=(G/\GG)\circ\GG$ where $G/\GG\in\tG_{n-k,1}$,
thus is a wheel (or a product of wheels) and $\GG\in\tG_{k,2}$.
Let us cap-off the boundary faces by some new auxiliary
configuration spaces. Recall that passing to $\ln(Z)$ instead of
$Z$, we can reduce our consideration only to connected graphs;
thus it would suffice to produce a $2n$-dimensional space $\hC_G$
with $\dd\hC_G=C_G$ for each wheel $G\in\tG_{n,1}$.

\begin{figure}[htb]
\centerline{\fig{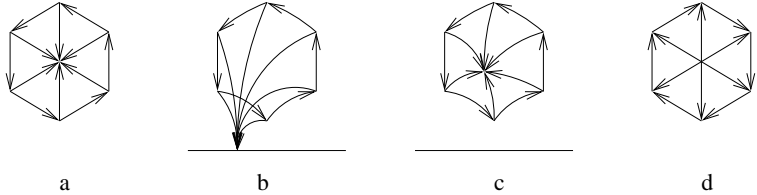,height=1in}}
\caption{Wheels and related configuration spaces}
\label{wheels_fig}
\end{figure}

For this purpose let us take a space of all embeddings
of the vertex set of $G$ in $\oH$, modulo the
2-dimensional group $Th(\oH)$ of translations and
homotheties of $\oH$, and construct its compactification
$\hC_G$ similarly to Section \ref{sub:gauss}.
This time there are no restrictions on the location
of the external vertex $v$ and it is free to run over $\H$.
Thus $\dim(\hC_G)=2(n+1)-2=2n$, and $C_G$ can be
realized as an infinity boundary face of $\hC_G$
corresponding to all embeddings for which $v\in\R$
and the other vertices remain in $\H$.
An example of a wheel $G\in\tG_{6,1}$ and points in
the corresponding configuration spaces $C_G$ and
$\hC_G$ are shown in Figures \ref{wheels_fig}a,
\ref{wheels_fig}b, and \ref{wheels_fig}c,
respectively.

If $\hC_G$ would have no other boundary faces, we
would be done.
So, let us try to get rid of all other strata of $\dd\hC_G$.
Other infinity faces are easy to get rid of: when some
internal vertices approach $\R$ (probably together with
$v$), some edges $e$ map under $\phi_e$ to $0$, so the
image of $\hat{\Phi}_G$ belongs to $O$ and the corresponding
face may be relativised.
Also, we can easily get rid of almost all internal
(principal and hidden) faces: when vertices of a proper
subgraph $\GG$ of $G$ collide together, $G/\GG$ contains
some multiple edges, so $G/\GG$ and hence also
$G/\GG\circ\GG$ vanish due to the antisymmetry
relation.
But when {\em all} vertices of $G$ collide together,
we get the problematic {\em anomaly} face $\hC_G(\R^2)$,
which is the space of all embeddings of the vertex
set of $G$ in $\R^2$ modulo the action of $Th(\R^2)$.

There are various ways to deal with the anomaly face.
One of them (not very elegant, but working) is as
follows: let us consider a graph $G'$ obtained from
$G$ by the reversal of directions of all edges, see
Figure \ref{wheels_fig}d.
The corresponding compactified configuration space 
$\hC_{G'}$ of all embeddings of the vertex set of $G'$ 
in $\H$ is of course the same $\hC_G$.
However, the map $\hat{\Phi}_{G'}$ differs from $\hat{\Phi}_G$.
This time we may relativise it on all faces (including
the infinite face where only $v$ approaches $\R$),
except for the anomaly face corresponding to collision
of all vertices. But on this anomaly face $\hat{\Phi}_{G'}$
coincides with $\hat{\Phi}_G$ up to reflections of some
copies of $S^1$ in $\T^{2n}$. So, folding each copy
of $S^1$ in two by identifying $z$ with $-z$ (or, 
equivalently, doubling the angle map $a\to 2a$), we
get the maps $2\Phi_{G'}$ and $2\Phi_G$ to agree on
the anomaly face.
Thus we may glue $\oC_{G'}$ and $\oC_G$ along the
anomaly face, getting a space with the only
non-relativised boundary face being the required one.
This finishes the cap-off construction.

Let me illustrate an application of this idea on a simple 
example of graphs with $n=2$ external vertices.

\subsection{Example: 2-spiked wheel}\label{sub:twospike}
\begin{figure}[htb]
\centerline{\fig{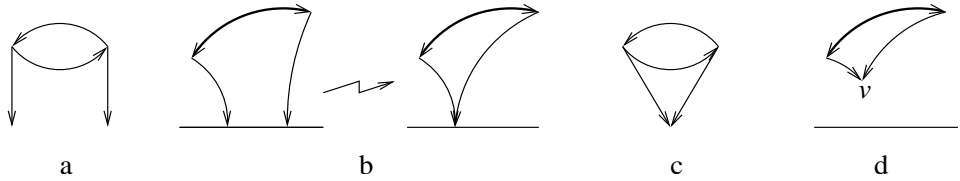,width=5in}}
\caption{Graphs with two internal vertices and their embeddings}
\label{two_fig}
\end{figure}

There is only one (up to a choice of an orientation) graph
$\GG\in\cG_{2,2}$ with $\ell(\GG)=1$, see Figure \ref{two_fig}a. The
only boundary face of $C^1_{2,2}$ which was not relativized or
zipped-up corresponds to collision of the two external vertices of
$\GG$, see Figures \ref{two_fig}b. This face may be identified with
$C_G$, where $G\in\cG_{2,1}$ is a wheel with two spikes, shown
in Figure \ref{two_fig}c. We cap-off this boundary face as 
described above, using an auxiliary 4-dimensional configuration 
space $\hC_G$, see Figure \ref{two_fig}d.  
The only boundary faces of $\hC_G$ which are not mapped to the
singular locus $O\in\T^4$ are $\C_G$ (correspondng to the vertex 
$v$ of $G$ approaching $\R$) and the anomaly face $\hC_G(\R^2)$  
(corresponding to all vertices of $G$ colliding together).
In this example we can get rid of the anomaly face without using 
$\hC_{G'}$, simply by folding each copy of $S^1$ in $T^4$ in two 
(identifying $\phi$ with $\phi+\pi$). 

After we get rid of the last boundary face $C_G$ of $C_\GG$ by 
gluing to it $\hC_G$, we may apply the scheme of Section 
\ref{sub:compare-good} to compute the relative degree of the 
map $\Phi_\GG\cup\hat{\Phi}_G:(C_\GG\cup\hC_G,D)\to(\T^4,O)$.
For any volume form $\Go$ on $\T^4$ we have
$$\int_{C_\GG}\Phi^*(\Go)+\int_{\hC_G}\hat{\Phi}^*(\Go)=
\deg(\Phi\cup\hat{\Phi})$$
We conclude that the weight $W_\GG(\Go)$ of $\GG$ and the weight
$\hat{W}_G(Go)=\int_{\hC_G}\hat{\Phi}^*(\Go)$ of the wheel $G$ are 
related by $W_\GG(\Go)+\hat{W}_G(\Go)=\deg(\Phi_\GG\cup\hat{\Phi}_G)$.
In particular, we can use for $\Go$ the standard uniform volume 
form on $\T^4$; this gives Kontsevich's weights $W^{Ko}_\GG$ 
for $W_\GG(\Go)$ and (following Shoikhet \cite{Sh}) $\hat{W}^{Ko}_G=0$,
so $W^{Ko}_\GG=\deg(\Phi_\GG\cup\hat{\Phi}_G)$.
On the other hand, we can compute $\deg(\Phi_\GG\cup\hat{\Phi}_G)$  
counting preimages of any regular value of $\Phi_\GG\cup\hat{\Phi}_G$; 
this gives a purely combinatorial method for calculation of $W^{Ko}_\GG$.

\end{document}